\documentclass[12pt]{elsarticle}
\usepackage{amssymb}
\usepackage{bbold}
\input ulem.sty
\input diagxy
\xyoption{curve}

\def\lax{{\rm Lax}_{\rm N}}
\def\cat{\mathbb{C}{\rm at}}
\def\pos{\mathbb{P}{\rm os}}
\def\loc{\mathbb{L}{\rm oc}}
\def\topos{\mathbb{T}{\rm opos}}
\def\spaces{\mathbb{T}{\rm op}}
\def\tcat{{\rm Cat}}
\def\tpos{{\rm Pos}}
\def\tloc{{\rm Loc}}
\def\ttopos{{\rm Topos}}
\def\tspaces{{\rm Top}}
\def\set{{\rm Sets}}
\def\tov{\hbox{$\to<150>$\hskip -.2in \raisebox{1.7pt}{\tiny$\bullet$} \hskip .1in}}
\def\id{{\rm id}}
\def\idv{\id^{\hbox{\tiny$\bullet$}}}
\def\dotv{\raisebox{1.5pt}{\tiny$\bullet$}}
\def\tinycirc{\raisebox{1pt}{\tiny$\circ$}}

\title{The Glueing Construction and Double Categories}
\vfill 
\author{Susan Niefield}
\ead{niefiels@union.edu}
\address{Department of Mathematics, Union College, Schenectady, NY 12308, USA}

\begin{document} 
\begin{frontmatter}

\begin{abstract} 
We introduce Artin-Wraith glueing and locally closed inclusions in double categories.  Examples include locales, toposes, topological spaces, categories, and posets.  With appropriate assumptions, we show that locally closed inclusions are exponentiable, and the exponentials are constructed via Artin-Wraith glueing.  Thus, we obtain a single theorem establishing the exponentiability of locally closed inclusions in these five cases.
\end{abstract}

\begin{keyword}
double category \sep Artin-Wraith glueing \sep lax colimits \sep locally closed inclusions \sep exponentiability 

\MSC[2010] 18D05  \sep 18A20 \sep 18A30 \sep 18B25 \sep 18B30 \sep 06D22

\end{keyword}
\end{frontmatter}

\newtheorem{lem}{Lemma}[section]
\newtheorem{prop}[lem]{Proposition}
\newtheorem{cor}[lem]{Corollary}
\newtheorem{thm}[lem]{Theorem}
\newproof{pf}{Proof}
\newdefinition{ex}[lem]{Example}
\newdefinition{exs}[lem]{Examples}

\section{Introduction}

Introduced in the 1960's by Ehresmann \cite{Ehres} and considered more recently by Par\'{e} et al \cite{DPP,GP99,GP04} and Shulman \cite{Shulman}, double categories provide a setting in which two types of morphisms are present in a single category.  Examples include the double categories: $\cat$ of small categories, functors, and profunctors; $\pos$ of posets, order-preserving maps, and order ideals; $\loc$ of locales, locale morphisms, and finite meet-preserving maps; $\spaces$ of topological spaces, continuous maps, and finite intersection-preserving maps on their open set lattices; and $\topos$ of toposes, geometric morphisms, and left exact functors.  

There are several exponentiability characterizations in which double categories implicitly arise.  In 1981, generalizing a  theorem  \cite{thesis} about topological spaces, Niefield \cite{inclusions} showed that the exponentiable inclusions of sublocales and subtoposes are precisely the locally closed ones.  To define the exponentials of open and closed (and hence, locally closed) sublocales, the Artin/Wraith glueing construction \cite{TT} was used to establish an equivalence between the 2-category $\tloc/{\cal O}(\mathbb 2)$ of locales over the Sierpinski locale and one whose objects are finite meet-preserving maps $m\colon X\to<125>X'$ on locales and morphisms $m\to<125>n$ are squares 
$$\bfig
\place(150,162)[\scriptstyle\ge] 
 \square<300,300>[X`Y`X'`Y';f`m`n`f']
\efig$$
where $f$ and  $f'$ are locale morphisms.  An analogous construction was established for the 2-category $\ttopos/{\cal S}^{\mathbb 2}$ of toposes over the Sierpinski topos.  In each case, both types of morphisms arise, and the squares are 2-cells in the related double categories $\loc$ and $\topos$.  In a 2001 unpublished note, Street used Benabou's equivalence  (see \cite{Street})
$$\tcat/B\simeq  \lax(B,{\cal P}rof)$$ where the latter is a 2-category whose objects are normal lax functors from a small category $B$ to the bicategory ${\cal P}rof$, to prove that $X\to<125> B$ is exponentiable in $\tcat$ (i.e., satisfies the Giraud/Conduch\'{e} condition \cite{Giraud,Conduche}) if and only if the corresponding $B\to<125> {\cal P}rof$ is a pseudofunctor.  This, too, can be seen as a result about the related double category $\cat$.  In fact, the locale/topos equivalences are analogous to Benabou's equivalence in the case where $B=\mathbb2$, since the morphisms of $\lax(B,{\cal P}rof)$ are functor-valued transformations which can be seen as squares of the form
$$\bfig
\place(150,155)[\to<80>] 
 \square<300,300>[X`Y`X'`Y';f`m`n`f']
\efig$$
where $f$ and $f'$ are functors and $m$ and $n$ are profunctors.  Note that, following \cite{TT}, our bicategory ${\cal P}rof$ is the dual of the bicategory ${\cal M}od$ considered in \cite{Street}.

We will see that these equivalences are a special case of a more general result.  For certain double categories $\mathbb D$ (including all five examples mentioned above) there is an adjoint pair $$\lax(B,\mathbb D)\two/->`<-/<150>^{\Gamma}_{\Delta}{\rm H}\mathbb D$$ where $\lax(B,\mathbb D)$ is a 2-category whose objects are normal lax functor from $B$ to $\mathbb D$ (see Section 4) and ${\rm H}\mathbb D$ is the horizontal category of $\mathbb D$.  Moreover, this adjunction factors through   ${\rm H}\mathbb D/\Gamma T$, when $\mathbb D$ has a (horizontal) terminal object $T$, and  gives an equivalence $$\lax(B,\mathbb D)\simeq {\rm H}\mathbb D/\Gamma T\eqno (\star)$$ when $B$ is a finite poset.  In particular, $\Gamma T=B$, when $\mathbb D =\cat$ and $\Gamma T$ is the Sierpinski locale (respectively, topos), when $B={\!\!\mathbb 2}$ and $\mathbb D =\loc$ (respectively, $\topos$).  Thus, we obtain a single theorem, not merely an analogy, which applies to the five double categories, and gives new results for $\loc$, $\spaces$, $\topos$, and $\pos$.

The common ground for the exponentiability results is that of inclusions of subobjects. In 1978, Niefield \cite{thesis,cart} showed that the inclusion of a subspace is exponentiable if and only if it is locally closed, and generalized this result to locales and toposes in \cite{inclusions}.  Also, the inclusion of a full subcategory (respectively, subposet) is exponentiable if and only if it is locally closed  \cite{BN} (respectively, \cite{posets}).  In each case, a subobject is called locally closed if it is the intersection of an open and a closed one, but open and closed are not categorical concepts.  We will define the notion of a locally closed inclusion in a double category $\mathbb D$, which agrees with the definitions in these five categories.  Using the equivalence $(\star)$ with $B=\mathbb 2$,  we will prove that locally closed inclusions are exponentiable in $ {\rm H}\mathbb D/D$, and thus obtain a construction of these exponentials via a single theorem.

We begin with an introduction to double categories (\S2) and the examples under consideration (\S3).  After a review of lax colimits in Section~4, we introduce a notion we call $B$-glueing (\S5) which captures the existence of the equivalence $(\star)$.  In Section~6, we define locally closed inclusions and establish their exponentiability in certain double categories (known as framed bicategories \cite{Shulman}).  We conclude (\S7) by showing that the equivalence $(\star)$ holds for all finite posets $B$, in the case where $\mathbb D$ is a framed bicategory with $\mathbb 2$-glueing and $B$-indexed lax colimits.

\bigskip
The author would like to thank Bob Par\'{e} and Geoff Cruttwell for pointing me in the direction of double categories.

\section{Double Categories} This section is a review of the properties of double categories used in this paper.  For details and more on the subject see Grandis/Par\'{e}  \cite{GP99,GP04} and Shulman \cite{Shulman}.

\medskip
A (pseudo) {\it double category} $\mathbb D$ is a pseudo category object in $\rm CAT$, i.e., functors $$\bfig 
\morphism(0,0)|a|/@{>}@<5pt>/<550,0>[\mathbb D_1\times_{\mathbb D_0}\mathbb D_1`\mathbb D_1;\pi_2]
\morphism(0,0)|m|/@{>}/<550,0>[\mathbb D_1\times_{\mathbb D_0}\mathbb D_1`\mathbb D_1;c]
\morphism(0,0)|b|/@{>}@<-5pt>/<550,0>[\mathbb D_1\times_{\mathbb D_0}\mathbb D_1`\mathbb D_1;\pi_1]
\morphism(550,0)|a|/@{>}@<5pt>/<350,0>[\mathbb D_1`\mathbb D_0;d_0]
\morphism(900,0)|m|<-350,0>[\mathbb D_0`\mathbb D_1;i]
\morphism(550,0)|b|/@{>}@<-5pt>/<350,0>[\mathbb D_1`\mathbb D_0;d_1]
\efig$$
such that $d_0 i=\id_{\mathbb D_0}, d_1 i=\id_{\mathbb D_0},d_0 c=d_0\pi_2,d_1 c=d_1\pi_1$, together with left and right identity isomorphisms ($\lambda$ and $\rho$, respectively) and an associativity isomorphism satisfying the usual coherence conditions.

The objects $X$ and morphisms $f\colon X\to<125>Y$ of $\mathbb D_0$ are called the {\it objects} and {\it horizontal morphisms} of $\mathbb D$.  The objects of $\mathbb D_1$ are called the {\it vertical morphisms} of $\mathbb D$, and denoted by $m\colon X\tov X'$, where $X=d_0m$ and $X'=d_1m$.  The morphisms $\varphi\colon m\to<125> n$ of $\mathbb D_1$ are called {\it cells} of $\mathbb D$, and denoted by $\varphi\colon m \to<125>^f_{f'}n$ or 
$$\bfig
\square<300,300>[X`Y`X'`Y';f`m`n`f']
\place(0,180)[\hbox{\tiny$\bullet$}]
\place(300,180)[\hbox{\tiny$\bullet$}]
\place(150,155)[\varphi]
\efig$$
Horizontal composition of morphisms (respectively, cells)  is given by composition in $\mathbb D_0$ (respectively, $\mathbb D_1$), and denoted by $\tinycirc$.  Vertical composition is given by $c$, and denoted by $\dotv$.  The cell compositions in 
$$\bfig
\iiixiii<300,300>[X`Y`Z`X'`Y'`Z'`X''`Y''`Z'';```````````]
\place(150,160)[\varphi']
\place(150,460)[\varphi]
\place(450,160)[\psi']
\place(450,460)[\psi]
\place(0,180)[\hbox{\tiny$\bullet$}]
\place(0,480)[\hbox{\tiny$\bullet$}]
\place(300,180)[\hbox{\tiny$\bullet$}]
\place(300,480)[\hbox{\tiny$\bullet$}]
\place(600,180)[\hbox{\tiny$\bullet$}]
\place(600,480)[\hbox{\tiny$\bullet$}]
\efig$$
are related by the {\it middle four interchange law} $$(\psi'\tinycirc\varphi')\dotv(\psi\tinycirc\varphi)=(\psi'\dotv\psi)\tinycirc(\varphi'\dotv\varphi)$$  The horizontal and vertical identity morphisms on $X$ are denoted by $\id_X$ and $\idv_X$, respectively.  A cell $\varphi$ is called {\it special} if both vertical or both horizontal morphisms are identities.

The objects, horizontal morphisms, and special cells 
$$\bfig
\square<300,300>[X`Y`X`Y;f`\idv_X`\idv_Y`g]
\place(0,180)[\hbox{\tiny$\bullet$}]
\place(300,180)[\hbox{\tiny$\bullet$}]
\place(150,155)[\varphi]
\efig$$ form a 2-category ${\rm H}\mathbb D$, called the {\it horizontal 2-category} of $\mathbb D$.  There is an analogous {\it vertical bicategory} ${\rm V}\mathbb D$. 

A {\it horizontal terminal object}  of $\mathbb D$ is an object $T$ such that there is a unique horizontal morphism $t\colon X\to<125>T$, for every object $X$, and a unique cell 
$$\bfig
\square<300,300>[X`T`X'`T;t`m`\idv_T`t']
\place(0,180)[\hbox{\tiny$\bullet$}]
\place(300,180)[\hbox{\tiny$\bullet$}]
\place(150,155)[\tau]
\efig$$ for every vertical morphism $m$.  A {\it horizontal initial object} is defined dually.

A {\it (vertical) companion} for $f\colon X\to<125>Y$ is a morphism $f_*\colon X\tov Y$ together with cells 
$$\bfig
\square(0,0)<300,300>[X`X`X`Y;\id_X`\idv_X`f_*`f]
\place(0,180)[\hbox{\tiny$\bullet$}]
\place(300,180)[\hbox{\tiny$\bullet$}]
\place(150,155)[\eta]
\square(1000,0)<300,300>[X`Y`Y`Y;f`f_*`\idv_Y`\id_Y]
\place(1000,180)[\hbox{\tiny$\bullet$}]
\place(1300,180)[\hbox{\tiny$\bullet$}]
\place(1150,155)[\varepsilon]
\efig$$
such that $\varepsilon\tinycirc \eta=\id_f$ and $\lambda_{f_*}\tinycirc(\varepsilon\dotv\eta)=\rho_{f_*}$, where  $\lambda$ and $\rho$ are the identity isomorphisms. 
A {\it (vertical) conjoint} for $f\colon X\to<125>Y$ is a morphism $f^*\colon Y\tov X$ together with cells 
$$\bfig
\square(0,0)<300,300>[X`Y`X`X;f`\id_X`f^*`\id_X]
\place(0,180)[\hbox{\tiny$\bullet$}]
\place(300,180)[\hbox{\tiny$\bullet$}]
\place(150,155)[\alpha]
\square(1000,0)<300,300>[Y`Y`X`Y;\id_Y`f^*`\idv_Y`f]
\place(1000,180)[\hbox{\tiny$\bullet$}]
\place(1300,180)[\hbox{\tiny$\bullet$}]
\place(1150,155)[\beta]
\efig$$
such that $\beta\tinycirc \alpha=\idv_f$ and $\rho_{f^*}\tinycirc(\alpha\dotv\beta)=\lambda_{f^*}$.  Moreover, $\mathbb D$ is called a {\it framed bicategory} (in the sense of  \cite{Shulman}) if and only if every horizontal morphism has a companion and conjoint. 

Note that if $f$ has a companion $f_*$ and a conjoint $f^*$, then $f_*$ is left adjoint to $f^*$ in the bicategory ${\rm V}\mathbb D$.   The notation $f_*$ (respectively, $f^*$) is justified by the fact that if $f$ has a companion (respectively, conjoint), then it is unique up to unique isocell.

\section{Examples} 
The following five double categories will be considered throughout this paper.  Each has horizontal initial  object, and every morphism has a companion and conjoint.  All but the last one, $\spaces$, have been previously considered.

\begin{ex} $\cat$ has small categories as objects, functors $f\colon X\to<125>Y$ and profunctors $m\colon X\tov X'$ (i.e., functors $m\colon X^{op} \times X'\to<125>\set$) as morphisms, and natural transformations $\varphi\colon f_*'\dotv m\to<125>n\dotv f_*$ as cells
$$\bfig
\square<300,300>[X`Y`X'`Y';f`m`n`f']
\place(0,180)[\hbox{\tiny$\bullet$}]
\place(300,180)[\hbox{\tiny$\bullet$}]
\place(150,155)[\varphi]
\efig$$
 The companion  $f_*\colon X\tov Y$ and conjoint  $f^*\colon Y\tov X$ of $f\colon X\to<125>Y$ are defined by $f_*(x,y)=Y(fx,y)$ and $f^*(y,x)=Y(y,fx)$, respectively.  Note that the transformation  $\varphi\colon f'_*\dotv m\to<125>n\dotv f_*$ can also be given by its mate  $\hat\varphi\colon m\dotv f^* \to<125>f'^*\dotv n$. The empty category $\mathbb 0$ and the one morphism category $\mathbb 1$ are the horizontal initial and terminal objects. 
\end{ex}

\begin{ex} $\pos$ has posets as objects, order-preserving maps $f\colon X\to<125>Y$ and order ideals $m\colon X\tov X'$ (i.e., upward closed subsets $m\subseteq X^{op}\times X'$) as morphisms.  Vertical composition is given by composition of relations.  There is a cell $\varphi\colon m \to<125>^f_{f'}n$ if and only if 
$$\bfig
\square<300,300>[X`Y`X'`Y';f`m`n`f']
\place(0,180)[\hbox{\tiny$\bullet$}]
\place(300,180)[\hbox{\tiny$\bullet$}]
\place(150,155)[\le]
\efig$$
or equivalently, $(x,x')\in m \Rightarrow (fx,f'x')\in n$.  The empty poset $\mathbb 0$ and the one element poset $\mathbb 1$ are the horizontal initial and terminal objects.  The companion  $f_*\colon X\tov Y$ and conjoint  $f^*\colon Y\tov X$  are given by $f_*=\{(x,y)\mid fx\le y\}$ and $f^*=\{(y,x)\mid y\le fx\}$, respectively. Note that there is a cell $\varphi\colon m\to<125>^f_{f'}n$ if and only if  $f'_*\dotv m\le n\dotv f_*$, or equivalently  $m\dotv f^*\le f'^*\dotv n$.
Also, $(\id_X)_*$  is the vertical identity $\idv_X$.
 \end{ex}

\begin{ex} $\loc$ has locales as objects, locale morphisms $f\colon X\to<125>Y$ and finite meet-preserving maps $m\colon X\tov X'$ as morphisms, and a cell $\varphi\colon m\to<125>^f_{f'}n$ if and only if 
$$\bfig
\square<300,300>[X`Y`X'`Y';f`m`n`f']
\place(0,180)[\hbox{\tiny$\bullet$}]
\place(300,180)[\hbox{\tiny$\bullet$}]
\place(150,155)[\ge]
\efig$$
i.e., $n\dotv f_*\le f'_*\dotv m$ or equivalently, $f'^*\dotv n\le m\dotv f^*$.  The initial locale $\mathbb 1$ and the terminal locale $\Omega$ are the horizontal initial and terminal objects.  The companion and conjoint of $f$ are the direct and inverse images $f_*$ and $f^*$, respectively. 
 \end{ex}

\begin{ex}
$\topos$ has (elementary) toposes $\cal X$ as objects, geometric morphisms $f\colon {\cal X}\to<125>{\cal Y}$ and  left exact functors $m\colon {\cal X}\tov{\cal X'}$ as morphisms, and natural transformations $\varphi\colon f'^*\dotv n\to<125>m\dotv f^*$, or equivalently, $\hat \varphi\colon n\dotv f_*\to<125>f'_*\dotv m$ as cells
$$\bfig
\square<300,300>[{\cal X}`{\cal Y}`{\cal X'}`{\cal Y'};f`m`n`f']
\place(0,180)[\hbox{\tiny$\bullet$}]
\place(300,180)[\hbox{\tiny$\bullet$}]
\place(150,160)[\varphi]
\efig$$
The one object topos is the horizontal initial  object, and the direct and inverse images of $f$ are the companion and conjoint, respectively.  To obtain a horizontal terminal object, it is necessary to restrict to Grothendieck toposes, or more generally, those bounded over a fixed base $\cal S$.
\end{ex}

\begin{ex} $\spaces$ has topological spaces $X$ as objects, continuous maps $f\colon X\to<125>Y$ as horizontal morphisms and finite intersection-preserving maps $m\colon{\cal O}(X)\to<125>{\cal O}(X')$ as vertical morphisms $m\colon X\tov X'$.  There is a cell
$$\bfig
\square<300,300>[X`Y`X'`Y';f`m`n`f']
\place(0,180)[\hbox{\tiny$\bullet$}]
\place(300,180)[\hbox{\tiny$\bullet$}]
\place(150,160)[\varphi]
\efig$$
if and only if there is a cell in $\loc$ of the form
$$\bfig
\square<400,300>[{\cal O}(X)`{\cal O}(Y)`{\cal O}(X')`{\cal O}(Y');f`m`n`f']
\place(0,180)[\hbox{\tiny$\bullet$}]
\place(400,180)[\hbox{\tiny$\bullet$}]
\place(200,155)[\ge]
\efig$$
The empty space and the one-point space are the horizontal  initial and terminal objects.  The companion and conjoint of $f\colon X\to<125>Y$ agree with those of  $f\colon {\cal O}(X)\to<125>{\cal O}(Y)$ in $\loc$.
 \end{ex}

\section{Lax Colimits}

This section begins with a review of the normal lax functors, transformations, and modifications that make up the 2-category $\lax(B,\mathbb D)$, for a small category $B$ and  a double category $\mathbb D$, and then  goes on to establish properties of $\lax(B,\mathbb D)$ which apply to the five examples under consideration.  

\medskip
A {\it vertical normal lax  functor} $F\colon B\to<125>\mathbb D$ consists of an object $Fb$, for each object $b$ of $B$, a vertical morphism $F_\beta\colon Fb\tov Fb'$, for every morphism  $\beta\colon b\to<125>b'$ of $B$, such that $F_{\id_b}=\idv_{Fb}$, and a special cell $\varphi_{\beta,\beta'}\colon F_{\beta'}\dotv F_\beta\to<125> F_{\beta'\beta}$ called a {\it comparison map}, for every composable pair,    satisfying the usual coherence conditions. 

A {\it horizontal lax transformation} $f\colon F\to<125> G\colon  B\to<125>\mathbb D$ consists of a horizontal morphism $f_b\colon Fb\to<125>Gb$, for every object $b$, and a cell
$$\bfig
\square<300,300>[Fb`Gb`Fb'`Gb';f_b`F_\beta`G_\beta`f_{b'}]
\place(0,180)[\hbox{\tiny$\bullet$}]
\place(300,180)[\hbox{\tiny$\bullet$}]
\place(150,160)[f_\beta]
\efig$$
for every  $\beta\colon b\to<125>b'$, such that $f_{\beta'\beta}\tinycirc \varphi_{\beta,\beta'}=\psi_{\beta,\beta'}\tinycirc(f_{\beta'}\dotv f_\beta)$, for all composable pairs, where $\varphi$ and $\psi$ are the comparison cells for $F$ and $G$, respectively.

A {\it modification} $\theta\colon f\to<125>g\colon F\to<125> G$ consists of a special cell  
$$\bfig
\square<300,300>[Fb`Gb`Fb`Gb;f_b`\idv_{Fb}`\idv_{Gb}`g_b]
\place(0,180)[\hbox{\tiny$\bullet$}]
\place(300,180)[\hbox{\tiny$\bullet$}]
\place(150,160)[\theta_b]
\efig$$
for every $b$, compatible with $f_\beta$ and $g_\beta$, for all $\beta\colon b\to<125>b'$.

Vertical normal lax functors $B\to<125>\mathbb D$, horizontal transformations, and modifications form a 2-category  denoted by $\lax(B,\mathbb D)$.

We say that $\mathbb D$ has {\it $B$-indexed (normal) lax colimits} if the constant functor $\Delta\colon{\rm H}\mathbb D\to<125>\lax(B,\mathbb D)$ has a left adjoint, which we denote by $\Gamma_{\!\!B}$ or just $\Gamma$.  One can show that $\Gamma F$ satisfies the universal property  
$$\bfig
\square/->`->`-->`->/<400,400>[Fb`\Gamma F`Fb'`Y;i_b`F_\beta`f`f_{b'}]
\place(0,250)[\hbox{\tiny$\bullet$}]
\place(220,340)[\scriptstyle i_\beta]
\morphism(0,0)<400,400>[Fb'`\Gamma F;]
\morphism(0,400)<400,-400>[Fb`Y;]
\place(220,80)[\scriptstyle f_\beta]
\efig$$
i.e., there is a family of horizontal morphisms $i_b\colon Fb\to<125>\Gamma F$  together with  cells $i_\beta\colon F_\beta\to<125>^{i_b}_{i_{b'}}\idv_{\Gamma F}$ which are compatible with the comparison maps, and universal in the following sense.  Given any object $Y$ and a family of horizontal morphisms $f_b\colon Fb\to<125>Y$ together with  cells $f_\beta\colon F_\beta\to<125>^{f_b}_{f_{b'}}\idv_Y$ which are compatible with the comparison maps, there exist a unique horizontal morphism  $f\colon \Gamma F\to<125>Y$ such that $f \tinycirc i_b =f_b$ and $\idv_f \tinycirc i_\beta =f_\beta$.  Furthermore, given any other such family $g_b\colon Fb\to<125>Y$ and cells $\theta_b\colon \idv_{Fb}{\to<125>^{f_b}_{g_b}}\idv_Y$ which are compatible with $f_\beta$ and $g_\beta$, there exists a unique cell $\theta\colon f\to<125>g$ such that $\theta\tinycirc \idv_{i_b}=\theta_b$.

If $\mathbb D$ has $B$-indexed lax colimits and $F\colon B\to<125>\mathbb D$ is a vertical normal lax functor, then $\Gamma$ induces a 2-functor $$\Gamma/F\colon \lax(B,\mathbb D)/F\to<125>{\rm H}\mathbb D/\Gamma F$$  Note that if $\mathbb D$ has a horizontal terminal object $T$, then $\lax(B,\mathbb D)$ has a terminal object $T$ defined by $Tb=T$ and $T_\beta=\idv_T$, and giving rise to a 2-functor $$\Gamma\colon \lax(B,\mathbb D)\to<125>{\rm H}\mathbb D/\Gamma T$$

\begin{ex}
$\cat$ has  all small lax colimits (also known as collages) given by the following Grothendieck construction.  For $F\colon B\to<125>\cat$, objects of $\Gamma F$ are pairs $(x,b)$, where $b$ is an object of $B$ and $x$ is an object of $Fb$.  Morphisms $(x,b)\to<125>(x',b')$ are pairs $(\alpha,\beta)$, where $\beta\colon b\to<125>b'$ and $\alpha\in F_{\beta}(x,x')$.  The identity morphism on $(x,b)$ is $(\id_x,\id_b)$, where the former is the identity on $x$ in $Fb$, and composition of morphisms is defined using the comparison maps.  The functors $i_b\colon Fb\to<125>\Gamma F$ are defined by $i_bx=(x,b)$ and $i_b(\alpha)=(\alpha,\id_b)$.  The cells $i_\beta\colon F_\beta\to<125>^{i_b}_{i_{b'}}\idv_{\Gamma F}$ are induced by the function
$$F_\beta(x,x')\times \Gamma F((x',b'),(x'',b''))\to<125> \Gamma F((x,b),(x'',b''))$$ defined by
$(\alpha,(\alpha',\beta'))\mapsto (\alpha'\alpha,\beta'\beta)$. One can show that the lax colimit  $\Gamma T$ is the category $B$, and thus we get a 2-functor $\Gamma\colon \lax(B,\cat)\to<125>\tcat/B$.

\end{ex}

\begin{ex} $\pos$ has all poset-indexed lax colimits which are defined similarly.  The elements of  $\Gamma F$ are pairs $(x,b)$, such that $b\in B$ and $x\in Fb$, with $(x,b)\le(x',b')$, if $b\le b'$ and $(x,x')\in F^b_{b'}$, where $F^b_{b'}\subseteq Fb\times Fb'$ is the order ideal associated with $b\le b'$.  Since $F^b_b$ is the order on the poset $Fb$ and $F^{b'}_{b''}\dotv F^b_{b'}\subseteq  F^b_{b''}$, for all $b\le b'\le b''$, it follows that $\Gamma F$  is a poset.  The maps $i_b\colon Fb\to<125>\Gamma F$ are defined by $i_bx=(x,b)$, and there are the required cells $i^b_{b'} \colon F^b_{b'}\to<125>^{i_b}_{i_{b'}}\idv_{\Gamma F}$, since $i_bx\le i_{b'}x'$, whenever $(x,x')\in F^b_{b'}$. The lax colimit  $\Gamma T$ is the poset $B$, and so we get  $\Gamma\colon \lax(B,\pos)\to<125>\tpos/B$.
\end{ex}

\begin{ex} $\loc$ has all poset-indexed lax colimits which are defined as follows.  Elements of  $\Gamma F$ are families $(x_b)_{b\in B}$,  where $x_b\in Fb$ and $x_{b'}\le F^b_{b'}x_b$, for all $b\le b'$.  Note that  $\Gamma F$ is a locale with the point-wise order, since $$\Gamma  F = \left(\sum_{b\in B} Fb\right)_{\!\! \textstyle g}$$ for the co-nucleus $$g\bigl((x_b)_b\bigr) = \left(\bigwedge_{b\le b'}F^b_{b'}x_b \right)_{\!\! b'}$$ The morphism $i_b\colon Fb\to<125>\Gamma F$ is the composite of the coproduct inclusion $Fb\to<125> \sum Fb$ and $g\colon\sum Fb\to<125> (\sum Fb)_g$.  The lax colimit  $\Gamma T$ is the locale $\downarrow\!\!Cl(B)$ of downward closed subsets of $B$, and so we get $$\Gamma\colon \lax(B,\loc)\to<125>\tloc/\!\!\downarrow\!\!Cl(B)$$
\end{ex}

\begin{ex} $\spaces$ has all poset-indexed lax colimits which are defined as follows.  The elements of  $\Gamma F$ are pairs $(x,b)$, where $b\in B$ and $x\in Fb$, with $i_b\colon Fb\to<125>\Gamma F$ given by $i_bx=(x,b)$.  A subset $U$ is open in $\Gamma F$, if $U_b$ is open in $Fb$, for all $b$, and $U_{b'}\subseteq F^b_{b'}(U_b)$, for all $b\le b'$, where $U_b=i_b^{-1}(U)$.  The  cells $i^b_{b'} \colon F^b_{b'}\to<125>^{i_b}_{i_{b'}}\idv_{\Gamma F}$ arise by definition of open set.  Moreover, the lax colimit  $\Gamma T$ is the Alexandrov space on $B$ (i.e., with the topology of downward closed sets), and so we get $\Gamma\colon \lax(B,\spaces)\to<125>\tspaces/B$.
\end{ex}

\begin{ex} $\topos$ has finite lax colimits, but the relevant adjoints in the definition are  pseudo adjoints, and this notion of lax colimit agrees with that of \cite{TT}. In particular, $\Gamma\Delta\cal S$ is the topos ${\cal S}^B$ of co-presheaves on $B$, and so we get $\Gamma/\Delta{\cal S}\colon \lax(B,\topos)/\Delta{\cal S}\to<125>\ttopos/{\cal S}^B$.
\end{ex}

\begin{prop} If $\mathbb D$ has $B$-indexed lax colimits and ${\rm H}\mathbb D$ has finite limits, then  $\Gamma/F\colon\lax(B,\mathbb D)/F\to<125>{\rm H}\mathbb D/\Gamma F$ has a right adjoint, for every $F\colon B\to<125>\mathbb D$.
\end{prop}

\begin{pf} This follows from the fact that a functor to a slice category has a right adjoint if and only if it composition with the forgetful functor does (see \cite{cart}, for example).  
\end{pf}

Proposition~4.6 applies to all of our examples except $\topos$ (which lacks a terminal object and some other finite limits).  In the case of $\cat$, this is B\'{e}nabou's equivalence mentioned above in the introduction.  For $\topos$, one can restrict to the subcategory $\topos_{\cal S}$ of Grothendieck toposes, or more generally, those bounded over a base $\cal S$, and obtain an adjoint pair $$\lax(B,\topos_{\cal S})/F \two/->`<-/<125>\ttopos_{\cal S}/\Gamma F$$ which we will see is an equivalence when $B$ is a finite poset.  In the general case, we know 
$\Gamma/F\colon \lax(\mathbb 2,\topos)/F\to<125>\ttopos/\Gamma F$ is always a pseudo equivalence \cite{inclusions}, and this, too, will be extended to all finite posets $B$. 

In each of our examples, there is a horizontal initial object $I$ which is vertically both initial and terminal, and every morphism $X\to<125>I$ is an isomorphism.  In addition, given $u\colon I\to<125>X$, the companion $u_*$ and conjoint $u^*$ are the unique vertical morphisms $I\tov X$ and $X\tov I$, respectively.  Such an object will be called a {\it zero object}, and denoted by $\mathbb O$.  At this point, we need not assume that $\mathbb D$ has companions and conjoints, just that the unique morphism $u\colon I\to<125>X$ does.

The following lemma will give an explicit description of the effect of the right adjoint to $\Gamma/F$ on objects, when $\mathbb D$ has a zero object.  

\begin{lem} If $\mathbb D$ has $B$-indexed lax colimits and a zero object $\mathbb O$, then the evaluation functor $(\ )_b\colon \lax(B,\mathbb D)\to<125> {\rm H}\mathbb D$ has a left adjoint, for all $b\in B$.
\end{lem} 

\begin{pf}  Define $L_b\colon {\rm H}\mathbb D \to<125> \lax(B,\mathbb D)$ by 
$(L_bX)_b= X$ and $(L_bX)_{b'}= \mathbb O$, for $b'\ne b$, with the obvious definition on morphisms.  Then one easily shows that $L_b$ is left adjoint to $(\ )_b$, as desired.
\end{pf}

\begin{prop} If $\mathbb D$ has $B$-indexed lax colimits and a zero object, then for every vertical normal lax functor $F\colon B\to<125>\mathbb D$, the following pullback exists in ${\rm H}\mathbb D$ 
$$\bfig
\square<325,300>[X_b`X`Fb`\Gamma F;``p`i_b]
\efig$$
and the right adjoint to $\Gamma/F\colon\lax(B,\mathbb D)/F\to<125>{\rm H}\mathbb D/\Gamma F$ takes $p\colon X\to<125> Fb$ to a lax functor given by  $b\mapsto X_b$. 
\end{prop}

\begin{pf} Consider the functors $${\rm H}\mathbb D/Fb\to<125>^{L_b}\lax(B,\mathbb D)/F\to<125>^{\Gamma/F}{\rm H}\mathbb D/\Gamma F$$ whose composite is the functor  $\Sigma_{i_b}$ defined by composition with $i_b\colon Fb\to<125>\Gamma F$.  Since $L_b$ and $\Gamma F$ have right adjoints, so does $\Sigma_{i_b}$, and the result follows.
\end{pf}

Lemma~4.7 can also be used to obtain a left adjoint $L_b$ to $(\ )_b$, for each object $b$, when the latter is considered as a functor $\lax(B,\mathbb D)/F\to<125> {\rm H}\mathbb D/Fb$.  This adjunction  will be used in Section~6, in the case where $B$ is the poset $\mathbb2$.  One can show that if  $\mathbb D$ is a framed bicategory, then $(\ )_b$ also has a right adjoint, for every $b$, but the absolute case does not imply the relative one, and so the latter must be established for each base $F$.  For simplicity, we only include the case $B=\mathbb2$, as that is that will be needed later.

\begin{prop} Suppose $\mathbb D$ is a framed bicategory and $l\colon D_0\tov D_1$.  Then $$(\ )_k\colon\lax(\mathbb2,\mathbb D)/l \to<125> {\rm H}\mathbb D/D_k$$ has a right adjoint, for $k=0,1$.
\end{prop}

\begin{pf}  Consider $R_k(Y\to<125>^{q_k}D_k)$ defined by
$$\bfig
\square<350,350>[D_0`D_0`D_1`D_1;\id`l`l`\id]
\place(0,575)[\hbox{\tiny$\bullet$}]
\place(350,575)[\hbox{\tiny$\bullet$}]
\place(175,525)[\psi_0]
\place(175,150)[\id_l]
\square(0,350)<350,350>[Y`D_0`D_0`D_0;q_0`{q_0}_*`\idv`]
\place(0,225)[\hbox{\tiny$\bullet$}]
\place(350,225)[\hbox{\tiny$\bullet$}]
\square(1200,0)<350,350>[D_1`D_1`Y`D_1;\id`q_1^*`\id`q_1]
\place(1200,575)[\hbox{\tiny$\bullet$}]
\place(1550,575)[\hbox{\tiny$\bullet$}]
\place(1375,155)[\psi_1]
\place(1375,525)[\id_l]
\square(1200,350)<350,350>[D_0`D_0`D_1`D_1;\id`l`l`]
\place(1200,225)[\hbox{\tiny$\bullet$}]
\place(1550,225)[\hbox{\tiny$\bullet$}]
\efig$$
for $k=0,1$, respectively.  Since $\mathbb D$ is framed, each cell $\varphi_k$ below
$$\bfig
\Atriangle(0,350)/<-`->`->/<175,225>[Y`X_0`D_0;f_0`q_0`]
\square<350,350>[X_0`D_0`X_1`D_1;p_0`m`l`p_1]
\place(0,220)[\hbox{\tiny$\bullet$}]
\place(350,220)[\hbox{\tiny$\bullet$}]
\place(175,175)[\varphi_0]
\square(1200,225)<350,350>[X_0`D_0`X_1`D_1;p_0`m`l`p_1]
\place(1200,440)[\hbox{\tiny$\bullet$}]
\place(1550,440)[\hbox{\tiny$\bullet$}]
\Vtriangle(1200,0)/->`->`<-/<175,225>[X_1`D_1`Y;`f_1`q_1]
\place(1375,400)[\varphi_1]
\efig$$
factors through a unique cell $\hat\varphi_k$ 
$$\bfig
\square(-350,0)/->`->``->/<350,700>[X_0`Y`X_1`D_1;f_0`m``p_1]
\place(-350,350)[\hbox{\tiny$\bullet$}]
\place(-200,350)[\hat\varphi_0]
\square<350,350>[D_0`D_0`D_1`D_1;\id`l`l`\id]
\place(0,575)[\hbox{\tiny$\bullet$}]
\place(350,575)[\hbox{\tiny$\bullet$}]
\place(175,525)[\psi_0]
\place(175,150)[\id_l]
\square(0,350)<350,350>[Y`D_0`D_0`D_0;q_0`{q_0}_*`\idv`]
\place(0,225)[\hbox{\tiny$\bullet$}]
\place(350,225)[\hbox{\tiny$\bullet$}]
\square(1050,0)/->`->``->/<350,700>[X_0`D_0`X_1`Y;p_0`m``f_1]
\place(1050,350)[\hbox{\tiny$\bullet$}]
\place(1200,350)[\hat\varphi_1]
\square(1400,0)<350,350>[D_1`D_1`Y`D_1;\id`q_1^*`\id`\id]
\place(1400,575)[\hbox{\tiny$\bullet$}]
\place(1750,575)[\hbox{\tiny$\bullet$}]
\place(1575,525)[\id_l]
\place(1575,150)[\psi_1]
\square(1400,350)<350,350>[D_0`D_0`D_1`D_1;\id`l`l`]
\place(1400,225)[\hbox{\tiny$\bullet$}]
\place(1750,225)[\hbox{\tiny$\bullet$}]
\efig$$
and it follows that $R_k$ is right adjoint to $(\ )_k$.
\end{pf}

\section{Glueing}

Suppose $\mathbb D$ is a framed bicategory with $B$-indexed colimits.  Then $\mathbb D$ has {\it $B$-glueing over $F\colon B\to<125>\mathbb D$} if $\Gamma/F\colon \lax(B,\mathbb D)/F\to<125>{\rm H}\mathbb D/\Gamma F$ is an equivalence of 2-categories with pseudo-inverse which takes $p\colon X\to<125>Fb$ to the lax functor over $F$ given by $$b\mapsto X_b, \quad \beta\mapsto(X_b{\hbox{$\to<150>^{(i_b)_*}$\hskip -.25in \raisebox{1.7pt}{\tiny$\bullet$} \hskip .1in}} \ X_\beta {\hbox{$\to<150>^{(i_{b'})^*}$\hskip -.25in \raisebox{1.7pt}{\tiny$\bullet$} \hskip .1in}} \ X_{b'})$$ where $X_b$ and $X_\beta$ are given by the following pullbacks in ${\rm H}\mathbb D$ 
$$\bfig
\square<350,300>[X_b`X`Fb`\Gamma_{\!\!B} F;i_b``p`]
\square(900,0)<350,300>[X_\beta`X`\Gamma_{\!\!\mathbb 2} F_\beta`\Gamma_{\!\!B} F;``p`]
\efig$$
We also say $\mathbb D$ has {\it $B$-glueing},  if $\mathbb D$ has $B$-glueing over $F$, for all $F$.

\begin{prop} If $\mathbb D$ has a horizontal terminal object $T$ and  $B$-glueing over $T$, then $\mathbb D$ has $B$-glueing.
\end{prop}

\begin{pf} This follows from the equivalence $$({\rm H}\mathbb D/\Gamma T)/(\Gamma F\to<100>\Gamma T)\simeq {\rm H}\mathbb D/\Gamma F$$
\end{pf}

We know $\cat$ has $B$-glueing, for all small categories $B$, and similarly, $\pos$ has $B$-glueing, for all posets $B$.  In  \cite{inclusions}, it is shown that $\topos$ and $\loc$ have $\mathbb 2$-glueing, but the definition must be taken in the pseudo sense in the case of $\topos$.  We will see that (with suitable assumptions on $\mathbb D$ applicable to our five examples),  if $\mathbb D$ has $\mathbb 2$-glueing, then  $\mathbb D$ has $B$-glueing, for all finite posets $B$.

\begin{prop} $\spaces$ has $\mathbb 2$-glueing.
\end{prop}

\begin{pf}  Consider $\Gamma\colon \lax(\mathbb 2,\spaces)\to<125>\tspaces/\mathbb 2$, where the space $\mathbb 2=\{0,1\}$ is the Sierpinski space with $\{0\}$ open.   Objects of $\lax(\mathbb 2,\spaces)$ are vertical morphisms $ X_0\tov X_1$ of $\spaces$, or equivalently, finite meet-preserving maps $m\colon{\cal O}(X_0)\to<125> {\cal O}(X_1)$.  We know $\Gamma(m)$ is the space $X_0+_m \!X_1$, whose points are pairs $(x,j)$, where $x\in X_j$ and $j\in \mathbb 2$, and $U\subseteq X_0+_m \!X_1$ is open, if each $U_j$ is open in $X_j$ and $U_1\subseteq m(U_0)$.  Now, $\Gamma$ is clearly full and faithful.  To see that it is essentially surjective, suppose $p\colon X\to<125>\mathbb 2$ is continuous, and let $X_0=p^{-1}(0)$, $X_1=p^{-1}(1)$, and $m\colon{\cal O}(X_0)\to<125> {\cal O}(X_1)$ be defined by $$m(U_0)=(U_0\cup X_1)^{\textstyle\circ}\cap X_1$$ where $(\ )^{\textstyle\circ}$ denotes the interior operator on ${\cal O}(X)$.  Then the induced map $f\colon X_0+_m \!X_1\to<125> X$ is clearly a continuous bijection.  To see that $f$ is an open map, suppose $U=U_0\cup U_1$, where each $U_j$ is open in $X_j$ and $U_1\subseteq m(U_0)$.  Then $U=(U_0\cup X_1)^{\textstyle\circ} \cap (X_0\cup U_1)$,  since 
\begin{eqnarray*}
(U_0\cup X_1)^{\textstyle\circ} \cap (X_0\cup U_1) &=& ((U_0\cup X_1)^{\textstyle\circ} \cap X_0)\cup  ((U_0\cup X_1)^{\textstyle\circ} \cap U_1)  \\ &=& U_0\cup U_1=U
\end{eqnarray*}
and both are clearly open in $X$.  Thus, $\Gamma$ is an equivalence, as desired.
\end{pf}

\section{Locally Closed Inclusions and Exponentiability} 

A horizontal morphism $i_0\colon D_0\to<125> D$ (respectively, $i_1\colon D_1\to<125> D$) is called an {\it open inclusion} (respectively, {\it closed inclusion}) if there is a lax colimit diagram in $\mathbb D$ of the form 
$$\bfig
\square<300,300>[D_0`D`D_1`D;i_0`l`\idv_D`i_1]
\place(0,180)[\hbox{\tiny$\bullet$}]
\place(300,180)[\hbox{\tiny$\bullet$}]
\place(150,160)[i_\beta]
\efig$$
Note that these ``inclusions" are not assumed to be monomorphisms, however, we will see that they are when $\mathbb D$ has a zero object and $\Gamma_{\!\!\mathbb 2}$ is faithful.  We will also see that ${\rm H}\mathbb D$ has pullbacks along these inclusions under suitable conditions, and so it will make sense to talk about a {\it locally closed inclusion} as the pullback of an open and a closed inclusion.

Open and closed inclusions in $\spaces$, $\loc$, and $\topos$ are inclusions of open and closed subobjects, respectively, by the glueing construction.  In $\cat$ and $\pos$, they are the discrete fibrations and opfibrations considered in \cite{BN} and \cite{posets}.
 
If $\mathbb D$ has a zero object, on easily shows that $D$ is the lax colimit of $\mathbb O\tov D$, and also of $D\tov \mathbb O$, and so $\id_D\colon D\to<125>D$ is both open and closed.

\begin{prop} If $\mathbb D$ has a zero object and $\Gamma_{\!\!\mathbb 2}$ exists and is faithful, then open and closed inclusions are monomorphisms.
\end{prop}

\begin{pf} Suppose $i_0\colon D_0\to<125> D$ is an open inclusion  and $i_0f=i_0g$, for some $f,g\colon X\to<125>D_0$, and consider
$$\bfig
\square(300,0)<300,300>[D_0`D_0`\mathbb O`D_1;\id_{D_0}```]
\place(300,180)[\hbox{\tiny$\bullet$}]
\place(600,180)[\hbox{\tiny$\bullet$}]
\place(0,180)[\hbox{\tiny$\bullet$}]
\place(140,300)[\two<150>^f_g]
\square/`->``->/<300,300>[X`D_0`\mathbb O`\mathbb O;```]
\efig$$
where the squares are the unique cells arising from $\mathbb O$.  Applying $\Gamma_{\!\!\mathbb 2}$, we get $X\two<125>^f_gD_0\to<125>^{i_0}D$.  Since $\Gamma_{\!\!\mathbb 2}$ is faithful, it follows that $f=g$, as desired.  Similarly, every closed inclusion is a monomorphism.
\end{pf}

Recall that an object $Y$ is {\it exponentiable} in a category $\cal D$, if $X\times Y$ exists, for all $X$, and the functor $- \times Y\colon {\cal D}\to<125>{\cal D}$ has a right adjoint, usually denoted by $( \ )^Y$.  The objects $Z^Y$ are called {\it exponentials}.  A morphism $q\colon Y\to<125>D$ is called {\it exponentiable} if it exponentiable in the slice category ${\cal D}/D$.  We also say $Y$ is {\it exponentiable over $D$}, in the latter case, and write $r^q\colon Z^Y\to<125>D$ for the exponential, given $r\colon Z\to<125>D$.

In each of the five categories under consideration, the inclusion of a subobject is exponentiable if and only if it is locally closed \cite{cart,inclusions,BN,posets}.
We will see that open, closed, and hence, locally closed inclusions are exponentiable in ${\rm H}\mathbb D$, when $\mathbb D$ is a framed bicategory with a zero object and $\mathbb2$-glueing (though we do not have a converse in this generality), but first a lemma.   
 
\begin{lem} Suppose $\mathbb D$ is a framed bicategory with a zero object $\mathbb0$.  Then  $l_0\colon D_0\tov \mathbb O$ and $l_1\colon\mathbb O\tov D_1$  are exponentiable over $l\colon D_0\tov D_1$ in $\lax(\mathbb 2,\mathbb D)$.
\end{lem}

\begin{pf} Consider
$$\bfig
\morphism(0,0)/@{->}@<5pt>/<700,0>[\lax(\mathbb2,\mathbb D)/l`{\rm H}\mathbb D/D_k;(\ )_k]
\morphism(0,0)|b|/<-/<700,0>[\lax(\mathbb2,\mathbb D)/l`{\rm H}\mathbb D/D_k;R_k]
\morphism(700,0)/@{->}@<5pt>/<700,0>[{\rm H}\mathbb D/D_k`\lax(\mathbb2,\mathbb D)/l;L_k]
\morphism(700,0)|b|/<-/<700,0>[{\rm H}\mathbb D/D_k`\lax(\mathbb2,\mathbb D)/l;(\ )_k]
\efig$$
\end{pf}
where $L_k\dashv (\ )_k\dashv R_k$ by Lemma~4.7 and Proposition~4.9, for $k=0,1$.  Then it is not difficult to show that $L_k\circ (\ )_k=-\times l_k$, and it follows that $l_k$ is exponentiable over $l$ in $\lax(\mathbb 2,\mathbb D)$.

\begin{thm} If $\mathbb D$ is a framed bicategory with $\mathbb 2$-glueing and a zero object, then locally closed inclusions are exponentiable in ${\rm H}\mathbb D$.
\end{thm}
\begin{pf} Suppose  $m\colon D_0\tov D_1$ and $D=\Gamma_{\!\!\mathbb 2}m$.  Then $\lax(\mathbb 2, \mathbb D)/m$ is equivalent to ${\rm H}\mathbb D/D$, by definition of  $\mathbb 2$-glueing.  Applying Lemma~6.2, we know that open and closed inclusions are exponentiable in ${\rm H}\mathbb D$.  Since products of exponentiables are exponentiable in any category, the desired result follows.
\end{pf}

Although the definition of locally closed inclusion involves only $\mathbb 2$-indexed lax colimits, there is an alternative using the poset $\mathbb 3=\{0,1,2\}$, namely, $i_1\colon D_1\to<125> D$ could be called locally closed if there is a vertical normal lax functor $F\colon \mathbb 3\to<125> \mathbb D$ such that  $\Gamma_{\!\!\mathbb 3}F=D$, $F1=D_1$, and $i_1\colon F1\to<125>\Gamma_{\!\!\mathbb 3}F$.  One can show that the two definitions are equivalent in the presence of $\mathbb 3$-glueing. In the next section, we will see that  $\mathbb D$ has $\mathbb 2$-glueing and $B$-indexed lax colimits for some finite poset $B$ (in particular, for $B=\mathbb 3$), then  if $\mathbb D$ has $B$-glueing.

\section{The Equivalence}

In this section, we prove a general theorem showing that $\mathbb 2$-glueing implies $B$-glueing, in the case where $B$ is a finite poset, and applies to the five double categories under consideration.

Suppose $B$ is a category which is the collage of a profunctor $k\colon B_0 \tov \mathbb 1$, where $\mathbb 1$ is the one morphism category with object $b_1$ and  $k\colon B_0 \tov \mathbb 1$  is given by $k(b,b_1)=B(b,b_1)$, for all $b$ in $B_0$.  If $\mathbb D$  has lax colimits indexed by $B_0$ and $B$, and $F\colon B\to<125> \mathbb D$, then there is an induced morphism $i_0\colon \Gamma_{\!\!B_0}F_0\to<125> \Gamma_{\!\!B}F$ such that the diagram 
$$\bfig
\qtriangle<400,300>[Fb`\Gamma_{\!\!B_0}F_0`\Gamma_{\!\!B}F;j_b`i_b`i_0]
\efig$$
commutes, where $F_0$ is the restriction of $F$ to $B_0$  and $j_b$ and $i_b$ are the colimit morphisms for $\Gamma_{\!\!B_0}F_0$ and $\Gamma_{\!\!B}F$, respectively, with analogous diagrams of cells.

In each of the examples, $i_0\colon \Gamma_{\!\!B_0}F_0\to<125>\Gamma_{\!\!B}F$  and $i_{b_1}\colon Fb_1\to<125>\Gamma_{\!\!B}F$ are  open and  closed inclusions, respectively.  In fact, this is true in general, when $B$ is  finite poset.

\begin{lem} Suppose $B=\Gamma_{\!\!\mathbb 2}(B_0 \tov \mathbb 1)$ is  a poset and $\mathbb D$ is a framed bicategory with $\mathbb 2$-glueing and lax colimits indexed by $B_0$ and $B$.  If $F\colon B\to<125>\mathbb D$ is a vertical normal lax functor, then  $i_0\colon \Gamma_{\!\!B_0}F_0\to<125>\Gamma_{\!\!B}F$  and $i_{b_1}\colon Fb_1\to<125>\Gamma_{\!\!B}F$ are  open and  closed inclusions via $$l\colon \Gamma_{\!\!B_0}F_0{\hbox{$\to<150>^{{i_0}_*}$\hskip -.25in \raisebox{1.7pt}{\tiny$\bullet$} \hskip .1in}}
\Gamma_{\!\!B}F{\hbox{$\to<150>^{i^*_{b_1}}$\hskip -.2in \raisebox{1.7pt}{\tiny$\bullet$} \hskip .1in}} Fb_1$$
\end{lem}

\begin{pf} Consider the diagram 
$$\bfig
\square/->`->``->/<400,400>[\Gamma_{\!\!B_0}F_0`\Gamma_{\!\!B} F`Fb_1`X;i_0`l``f_{b_1}]
\place(0,250)[\hbox{\tiny$\bullet$}]
\place(220,330)[\scriptstyle \varphi]
\place(360,270)[\scriptstyle i_{b_1}]
\morphism(0,0)<400,400>[Fb_1`\Gamma_{\!\!B} F;]
\morphism(0,400)<400,-400>[\Gamma_{\!\!B_0}F_0`X;]
\place(220,85)[\scriptstyle \psi]
\place(360,135)[\scriptstyle f_0]
\efig$$ where $\varphi\colon l\to<125>^{i_0}_{i_{b_1}}\idv_{\Gamma_{\!\!B}F}$ corresponds to $\id_l\colon l\to<125>^{i_0}_{i_{b_1}}l$.  Note that, since $\mathbb D$ has $\mathbb 2$-glueing,  we can assume $F^b_{b'}=i^*_{b'}{i_b}_*$, for all $b\le b'$ in $B$.  Composing the diagram with the morphisms $j_b\colon Fb\to<125>\Gamma_{\!\!B_0}F_0$, for each $b\in B_0$, gives rise to a diagram
$$\bfig
\square/->`->`-->`->/<400,400>[Fb`\Gamma_{\!\!B} F`Fb'`X;i_b`F^b_{b'}`f`f_{b'}]
\place(0,250)[\hbox{\tiny$\bullet$}]
\place(220,330)[\scriptstyle \varphi^b_{b'}]
\morphism(0,0)<400,400>[Fb'`\Gamma_{\!\!B} F;]
\morphism(0,400)<400,-400>[Fb`X;]
\place(220,85)[\scriptstyle \psi^b_{b'}]
\efig$$
and hence, a unique morphism $f\colon \Gamma_{\!\!B} F\to<125>X$ such that $f\tinycirc i_b=f_b$, for all $b\in B$,  and $(\idv_f)\dotv \varphi^b_{b'}=\psi^b_{b'}$, for all $b\le b'$, and the desired result follows.
\end{pf}

\begin{thm} If $\mathbb D$ is a framed bicategory with $\mathbb 2$-glueing and lax colimits indexed by finite posets, then $\mathbb D$ has $B$-glueing, for all finite posets $B$.
\end{thm}

\begin{pf} We proceed by induction on $n$ to show that if $\mathbb D$ has $B$-glueing for all finite $n$-element posets $B$.  For the induction step, assume $|B|\ge2$, and let $B_1=\{b_1\}$, for some maximal element of $B$, and $B_0=B\setminus B_1$.  Then $B=\Gamma_{\!\!\mathbb2}(k)$, for $k\colon B_0\tov B_1$ defined by $k=\{(b,b_1)\mid b\le b_1 \ {\rm in} \ B\}$.  Given $F\colon B\to<125>\mathbb D$, consider $F_0=F|_{B_0}$, $F_1=F|_{B_1}$, and $l\colon \Gamma_{\!\!B_0}F_0\tov \Gamma_{\!\!B_1}F_1$, as defined in Lemma~7.1 above.  Then, applying the lemma, we see that the functor $\Gamma_B/F$
factors as 
$$\lax(B,\mathbb D)/F\simeq\lax(\mathbb 2,\mathbb D)/ l\simeq {\rm H}\mathbb D/\Gamma_{\!\!B}F$$ via equivalences with pseudo-inverses having the desired properties.
\end{pf}

Since $\spaces$ is such a framed bicategory with a terminal object $T$ and $\Gamma_{\!\!B}T=B$ (c.f. Example~4.4) and every finite $T_0$-space arises in this way, we get the following corollary.  

\begin{cor} If $B$ is a finite $T_0$-space, then $\lax(B,\spaces)\simeq \tspaces/B$,  where $B$ is a poset via the specialization order.
\end{cor}

Applying the theorem to $\loc$, where $\Gamma_{\!\!B}T=\downarrow\!\!Cl(B)$, yields:

\begin{cor} If $B$ is a finite poset, then $\lax(B,\loc)\simeq \tloc/\!\!\downarrow\!\!Cl(B)$.
\end{cor}

Finally, if $\cal S$ is a topos, we can apply the pseudo version of the theorem.
We know $\Gamma_{\!\!B}\Delta{\cal S}\simeq /{\cal S}^B$, and so:

\begin{cor} If $B$ is a finite poset, then  $\lax(B,\topos)/\Delta{\cal S}\simeq\ttopos/{\cal S}^B$, and  restricting to Grothendieck toposes, $\lax(B,\topos_{\cal S})\simeq\ttopos/{\cal S}^B$.
\end{cor}

\end{document}